\documentclass[a4paper]{amsart}

\usepackage[german,english]{babel}

\RequirePackage{amsmath}
\RequirePackage{bm}
\RequirePackage{amssymb}
\RequirePackage{upref}
\RequirePackage{amsthm}
\RequirePackage{enumerate}
\RequirePackage{pb-diagram}
\RequirePackage{amsfonts}
\RequirePackage[mathscr]{eucal}
\RequirePackage{verbatim}
\RequirePackage{xr}
\RequirePackage{graphicx}
\RequirePackage{floatflt}
\usepackage{calc}
\usepackage{xspace}

\newcommand{\cf}{cf.\@\xspace}
\newcommand{\resp}{resp.\@\xspace}

\newcommand{\al}{\alpha}
\newcommand{\bet}{\beta}
\newcommand{\ga}{\gamma}
\newcommand{\de}{\delta }
\newcommand{\e}{\epsilon}

\newcommand{\f}{\varphi}
\newcommand{\h}{\eta}

\newcommand{\ka}{\kappa}
\newcommand{\lam}{\lambda}

\newcommand{\n}{\nu}

\newcommand{\s}{\sigma}
\newcommand{\x}{\xi}

\newcommand{\C}{\varGamma}

\newcommand{\F}{\varPhi}

\newcommand{\Om}{\varOmega}

\newcommand{\di}[1]{#1\nobreakdash-\hspace{0pt}dimensional}%\di n
\newcommand{\nbdd}{\nobreakdash--}

\newcommand{\fm}[1]{F_{|_{M_#1}}}
\newcommand{\fv}[2]{#1\hspace{0pt}_{|_{#2}}}

\newcommand{\so}{{\mc S_0}}

\newcommand{\const}{\tup{const}}
\newcommand{\ndash}{\nobreakdash--}

\newcommand{\msp[1]}[1]{\mspace{#1mu}}
\newcommand{\low}[1]{{\hbox{}_{#1}}}

\newcommand{\R}[1][n+1]{{\protect\mathbb R}^{#1}}
\newcommand{\Hh}[1][n+1]{{\protect\mathbb H}^{#1}}

\newcommand{\N}{{\protect\mathbb N}}

\newcommand{\eR}{\stackrel{\lower1ex \hbox{\rule{6.5pt}{0.5pt}}}{\msp[3]\R[]}}
\newcommand{\eN}{\stackrel{\lower1ex \hbox{\rule{6.5pt}{0.5pt}}}{\msp[1]\N}}
\newcommand{\eO}{\stackrel{\lower1ex
\hbox{\rule{6pt}{0.5pt}}}{\msc O}}

\DeclareMathOperator{\graph}{graph}

\DeclareMathOperator{\id}{id}

\DeclareMathOperator{\grad}{grad}

\newcommand\ra{\rightarrow}

\newcommand\pa{\partial}
\newcommand\pde[2]{\frac {\partial#1}{\partial#2}}
\newcommand{\un}{\infty}
\newcommand{\A}{\forall}

\newcommand{\set}[2]{\{\,#1\colon #2\,\}}
\newcommand{\uu}{\cup}
\newcommand{\ii}{\cap}
\newcommand{\uuu}{\bigcup}

\newcommand{\uud}{ \stackrel{\lower 1ex \hbox {.}}{\uu}}
\newcommand{\uuud}[1]{ \stackrel{\lower 1ex \hbox {.}}{\uuu_{#1}}}
\newcommand\su{\subset}

\newcommand{\sminus}[1][28]{\raise 0.#1ex\hbox{$\scriptstyle\setminus$}}
\newcommand{\cpl}{\complement}

\newcommand{\wed}{\wedge}

\newcommand\ti{\times }

\newcommand{\abs}[1]{\lvert#1\rvert}

\newcommand{\norm}[1]{\lVert#1\rVert}

\newcommand{\spd}[2]{\protect\langle #1,#2\protect\rangle}

\newcommand\ch[3]{\varGamma_{#1#2}^#3}
\newcommand\cha[3]{{\bar\varGamma}_{#1#2}^#3}

\newcommand{\riem}[4]{R_{#1#2#3#4}}
\newcommand{\riema}[4]{{\bar R}_{#1#2#3#4}}

\newcommand{\tit}{\textit}

\newcommand{\tup}{\textup}% text upright

\newcommand{\mc}{\protect\mathcal}
\newcommand{\msc}{\protect\mathscr}

\providecommand{\bysame}{\makebox[3em]{\hrulefill}\thinspace}

\newcommand{\ci}{\cite}

\newcommand{\cq}[1]{\glqq{#1}\grqq\,}

\newcommand{\bt}{\begin{thm}}
\newcommand{\bl}{\begin{lem}}
\newcommand{\bc}{\begin{cor}}
\newcommand{\bd}{\begin{definition}}
\newcommand{\bpp}{\begin{prop}}
\newcommand{\br}{\begin{rem}}
\newcommand{\bn}{\begin{note}}
\newcommand{\be}{\begin{ex}}
\newcommand{\bes}{\begin{exs}}
\newcommand{\bb}{\begin{example}}
\newcommand{\bbs}{\begin{examples}}
\newcommand{\ba}{\begin{axiom}}
\newcommand{\bas}{\begin{assumption}}

\newcommand{\et}{\end{thm}}
\newcommand{\el}{\end{lem}}
\newcommand{\ec}{\end{cor}}
\newcommand{\ed}{\end{definition}}
\newcommand{\epp}{\end{prop}}
\newcommand{\er}{\end{rem}}
\newcommand{\en}{\end{note}}
\newcommand{\ee}{\end{ex}}
\newcommand{\ees}{\end{exs}}
\newcommand{\eb}{\end{example}}
\newcommand{\ebs}{\end{examples}}
\newcommand{\ea}{\end{axiom}}
\newcommand{\eas}{\end{assumption}}

\newcommand{\bp}{\begin{proof}}
\newcommand{\ep}{\end{proof}}
\newcommand{\eps}{\renewcommand{\qed}{}\end{proof}}

\newcommand{\bal}{\begin{align}}

\newcommand{\bi}[1][1.]{\begin{enumerate}[\upshape #1]}
\newcommand{\bia}[1][(1)]{\begin{enumerate}[\upshape #1]}
\newcommand{\bin}[1][1]{\begin{enumerate}[\upshape\bfseries #1]}
\newcommand{\bir}[1][(i)]{\begin{enumerate}[\upshape #1]}
\newcommand{\bic}[1][(i)]{\begin{enumerate}[\upshape\hspace{2\cma}#1]}
\newcommand{\bis}[2][1.]{\begin{enumerate}[\upshape\hspace{#2\parindent}#1]}
\newcommand{\ei}{\end{enumerate}}

\newcommand\ndots{\raise 0.47ex \hbox {,}\hskip0.06em\cdots %
     \raise 0.47ex \hbox {,}\hskip0.06em} 

\newcommand{\q}{\quad}
\newcommand{\qq}{\qquad}

\newcommand{\hp}{\hphantom}

\newcommand\nd{\noindent}

\newskip\Csmallskipamount                                                
\Csmallskipamount=\smallskipamount
\newskip\Cmedskipamount
\Cmedskipamount=\medskipamount
\newskip\Cbigskipamount
\Cbigskipamount=\bigskipamount

\newcommand\cvs{\vspace\Csmallskipamount}   
\newcommand\cvm{\vspace\Cmedskipamount}

\newskip\csa
\csa=\smallskipamount

\newskip\cma
\cma=\medskipamount

\newskip\cba
\cba=\bigskipamount

\newdimen\spt
\newcommand\citem{\cvs\advance\itemno by
1{(\romannumeral\the\itemno})\hskip3pt}
\newcommand{\bitem}{\cvm\nd\advance\itemno by
1{\bf\the\itemno}\hspace{\cma}}

\newcount\itemno
\newcommand{\las}[1]{\label{S:#1}}

\newcommand{\lae}[1]{\label{E:#1}}
\newcommand{\lat}[1]{\label{T:#1}}
\newcommand{\lal}[1]{\label{L:#1}}
\newcommand{\lad}[1]{\label{D:#1}}

\newcommand{\lar}[1]{\label{R:#1}}

\newcommand{\rs}[1]{Section~\ref{S:#1}}

\newcommand{\rt}[1]{Theorem~\ref{T:#1}}
\newcommand{\rl}[1]{Lemma~\ref{L:#1}}
\newcommand{\rd}[1]{Definition~\ref{D:#1}}

\newcommand{\rr}[1]{Remark~\ref{R:#1}}

\newcommand{\re}[1]{\eqref{E:#1}}

\newskip\thmskip
\thmskip=\parindent

\newskip\hsk
\newenvironment{hinw}{\labelsep=0pt\begin{list}{}{\labelsep=0pt\itemindent=0pt\labelwidth=0pt\leftmargin=\parindent\rightmargin=0pt\partopsep=\cba}%
\item\it\nopagebreak\nopagebreak}%
{\end{list}}

\newcommand\bh{\begin{hinw}}
\newcommand{\eh}{\end{hinw}}

\newtheoremstyle{normal}% name
  {\cba}%      Space above, empty = `usual value'
  {\cba}%      Space below
  {}% Body font
  {\thmskip}%Indent amount (empty = no indent, \parindent = para indent)
  {\bfseries}% Thm head font
  {.}%        Punctuation after thm head
  {\hsk}%     Space after thm head: " " = normal interword space;
  {}% Thm head spec

\newtheoremstyle{abschnitt}% name
  {\cba}%      Space above, empty = `usual value'
  {\cba}%      Space below
  {}% Body font
  {\thmskip}% Indent amount (empty = no indent, \parindent = para indent)
  {\bfseries}% Thm head font
  {.}%        Punctuation after thm head
  {\hsk}%     Space after thm head: " " = normal interword space;
  {}% Thm head spec

\newtheoremstyle{italic}% name
  {\cba}%      Space above, empty = `usual value'
  {\cba}%      Space below
  {\itshape}% Body font
  {\thmskip}%  Indent amount (empty = no indent, \parindent = para indent)
  {\bfseries}% Thm head font
  {.}%        Punctuation after thm head
  {\hsk}%     Space after thm head: " " = normal interword space;
  {}% Thm head spec

\newtheoremstyle{aufgaben}% name
  {\cba}%      Space above, empty = `usual value'
  {\cba}%      Space below
  {}% Body font
  {}%         Indent amount (empty = no indent, \parindent = para indent)
  {\normalsize\bfseries}% Thm head font
  {.}%        Punctuation after thm head
  {\hsk}%     Space after thm head: " " = normal interword space;
  {}% Thm head spec

\newtheoremstyle{break}% name
  {\cba}%      Space above, empty = `usual value'
  {\cba}%      Space below
  {\itshape}% Body font
  {}%         Indent amount (empty = no indent, \parindent = para indent)
  {\bfseries}% Thm head font
  {.}%        Punctuation after thm head
  {\newline}% Space after thm head: \newline = linebreak
  {}%         Thm head spec

\swapnumbers
\theoremstyle{italic}
\newtheorem{thm}[subsection]{Theorem}
\newtheorem{lem}[subsection]{Lemma}
\newtheorem{prop}[subsection]{Proposition}
\newtheorem{cor}[subsection]{Corollary}

\theoremstyle{normal}
\newtheorem{rem}[subsection]{Remark}
\newtheorem{definition}[subsection]{Definition}
\newtheorem{example}[subsection]{Example}
\newtheorem{examples}[subsection]{Examples}
\newtheorem{ex}[subsection]{Exercise}
\newtheorem{note}[subsection]{}
\newtheorem{axiom}[subsection]{Axiom}
\newtheorem{assumption}[subsection]{Assumption}

\theoremstyle{aufgaben}
\newtheorem{exs}[subsection]{Exercises}

\swapnumbers

\numberwithin{equation}{section}
\numberwithin{figure}{section}

\newenvironment{textequation}[1][0.8]
{\begin{equation}
\begin{aligned}
\begin{minipage}{#1\linewidth}}
{\end{minipage}
\end{aligned}
\end{equation}
\ignorespacesafterend}

\newcommand{\btext}{\begin{textequation}}
\newcommand{\etext}{\end{textequation}}

\usepackage{calc}

\newlength{\oddsidemarginlength}
\newlength{\topmarginlength}

\newcounter{numberoflines}
\newcounter{tempcc}
\oddsidemargin=\oddsidemarginlength
\evensidemargin=\oddsidemargin
\topmargin=\topmarginlength
\usepackage{hyperref}

\begin{document}

\flushbottom

%\larger[1]
%\frontmatter

\title[Minkowski type problems in hyperbolic space]{Minkowski type problems for convex hypersurfaces in hyperbolic space}

% author one information
\author{Claus Gerhardt}
\address{Ruprecht-Karls-Universit\"at, Institut f\"ur Angewandte Mathematik,
Im Neuenheimer Feld 294, 69120 Heidelberg, Germany}
%\curraddr{}
\email{gerhardt@math.uni-heidelberg.de}
\urladdr{http://www.math.uni-heidelberg.de/studinfo/gerhardt/}
\thanks{}

% author two information
%\author{}
%\address{}
%\curraddr{}
%\email{}
%\thanks{}
%
\date{\today}
%
% at present the "communicated by" line appears only in ERA and PROC
%\commby{}

%\dedicatory{}

\begin{abstract}
We consider the problem $F=f(\nu)$ for strictly convex, closed hypersurfaces in $\Hh$ and solve it for curvature functions $F$ the inverses of which are of class $(K^*)$.
\end{abstract}

\maketitle

\tableofcontents

\setcounter{section}{-1}
\section{Introduction}

In the classical Minkowski problem in $\R$ one wants to find a strictly convex closed hypersurface $M\su \R$ such that its Gau{\ss} curvature $K$ equals a given function $f$ defined in the normal space of $M$ or equivalently defined on $S^n$
\begin{equation}
\fv KM=f(\nu).
\end{equation}

The problem has been partially solved by Minkowski \cite{minkowski:minkowski}, Alexandrov \cite{alexandrov:minkowski}, Lewy \cite{lewy:minkowski}, Nirenberg \cite{nirenberg:minkowski}, and Pogorelov \cite{pogorelov:minkowski}, and in full generality by Cheng and Yau \cite{cheng-yau:minkowski}.

Instead of prescribing the Gaussian curvature other curvature functions $F$ can be considered, i.e., one studies the problem
\begin{equation}\lae{0.2}
\fv FM=f(\nu).
\end{equation}

If $F$ is one of the symmetric polynomials $H_k$, $1\le k\le n$, this problem has been solved by Guan and Guan \cite{guan:annals}.  They proved that \re{0.2} has a solution, if $f$ is invariant with respect to a fixed point free group of isometries of $S^n$. 

In a previous work \cite{cg:minkowski} we solved the problem \re{0.2} for strictly convex hypersurfaces $M\su S^{n+1}$ and for curvature functions $F$ the inverses of which are of class $(K)$, see \cite[Definition 1.3]{cg:indiana}. These $F$ include all $H_k$, $1\le k\le n$, $\abs{A}^2$, and also any symmetric, convex curvature function homogeneous of degree $1$, \cf \cite[Lemma 1.6]{cg97}.

In the present paper we consider the problem \re{0.2} for strictly convex hypersurfaces in $\Hh$ and for curvature functions $F$ the inverses of which belong to a subclass of $(K)$, the so-called class $(K^*)$, \cf \cite[Definition 1.6]{cg:indiana}.

Among the curvature functions $F$ that satisfy this requirement are the Gausssian curvature $F=K=H_n$, and all curvature functions that can be written as
\begin{equation}
F=H_k K^a,\qq1\le k\le n,
\end{equation}
where $a>0$ is a constant, as well as positive powers of those functions.

The Minkowski space $\R[n+1,1]$ contains two spaces of constant curvature as hypersurfaces, namely, $\Hh$ which is defined as
\begin{equation}
\Hh=\set{x\in\R[n+1,1]}{\spd xx=-1,\, x^0>0}
\end{equation}
and  the de Sitter space-time $N$, a Lorentzian manifold  of constant curvature $K_N=1$
\begin{equation}
N=\set{x\in\R[n+1,1]}{\spd xx=1}.
\end{equation}

We shall show in \rs{4} that for any closed strictly convex hypersurface $M\su \Hh$ there exists a Gau{\ss} map
\begin{equation}\lae{0.5b}
x\in M\ra \tilde x\in M^*\su N,
\end{equation}
where $M^*$ is the polar set of $M$. $M^*$ is spacelike, also strictly convex, as smooth as $M$, and the Gau{\ss} map is a diffeomorphism.

On the other hand, for any given closed, spacelike, connected,  strictly convex hypersurface $M\su N$ there also exists a Gau{\ss} map  
\begin{equation}
x\in M\ra \tilde x\in M^*\su \Hh
\end{equation}
which maps $M$ onto a closed, strictly convex hypersurface in hyperbolic space. These Gau{\ss} maps are inverse to each other.

If we consider $M\su \Hh$ as an embedding in $\R[n+1,1]$ of codimension $2$, so that the tangent spaces $T_x(M)$ and $T_x(\Hh)$ can be identified with  subspaces of $T_x(\R[n+1,1])$, then the image of  the point $x$ under the Gau{\ss} map is exactly the normal vector $\nu\in T_x(\Hh)$
\begin{equation}
\tilde x=\nu\in T_x(\Hh)\su T_x(\R[n+1,1]).
\end{equation}

Thus, the equation \re{0.2} can also be written in the form
\begin{equation}\lae{0.5}
\fv FM=f(\tilde x)\qq\A\,Ê x\in M,
\end{equation}
where $f$ is given as a function defined in $N$.

Using \re{0.5b} we shall prove that \re{0.5} has a dual problem, namely,
\begin{equation}
\fv {\tilde F}{M^*}=f^{-1}(\tilde x)\qq\A\,  \tilde x\in M^*,
\end{equation}
where $\tilde F$ is the inverse of $F$
\begin{equation}
\tilde F(\ka_i)=\frac1{F(\ka_i^{-1})}.
\end{equation}

In the dual problem the curvature is not prescribed by a function defined in the normal space, but by a function defined on the hypersurface.

Both problems are equivalent, solving one also leads to a solution of the dual one; notice also that
\begin{equation}
M^{**}=M\q\wed\q \tilde{\tilde x}=x.
\end{equation}

To find a solution  we assume that $(\tilde F,f^{-1})$ satisfy barrier conditions, \cf \rd{0.1} for details.

Then we shall prove 

\bt\lat{0.4}
Let $F\in C^{m,\al}(\C_+)$, $2\le m$, $0<\al<1$,  be a symmetric, positively homogeneous and monotone curvature function such that its inverse $\tilde F$ is of class $(K^*)$, let $0<f\in C^{m,\al}(N)$ and assume that the barrier conditions for $(\tilde F, f^{-1})$ are satisfied,  then the dual problems
\begin{equation}
\fv FM=f(\tilde x)
\end{equation}
and
\begin{equation}
\fv{\tilde F}{M^*}=f^{-1}(\tilde x)
\end{equation}
have strictly convex solutions $M$ \resp $M^*$ of class $C^{m+2,\al}$,  where $M^*$ is spacelike. 
\et

 The paper is organized as follows: \rs{1} gives an overview of the definitions and conventions we rely on.  In \rs{2} we define the Beltrami map from $\Hh$ to $\R$ with the help of which we prove Hadamard's theorem for strictly convex hypersurfaces in $\Hh$ in \rs{3}. 
 
The Gau{\ss} maps and their properties are treated in \rs{4}. In the last three sections we prove the existence of a solution in $N$ using a curvature flow method.  
\section{Notations and definitions}\las 1

The main objective of this section is to state the equations of Gau{\ss}, Codazzi,
and Weingarten for hypersurfaces. Since we are dealing with hypersurfaces in a Riemannian space as well as in a Lorenztian space, we shall formulate the governing equations of a
hypersurface $M$ in a semi-riemannian \di{(n+1)} manifold $N$, which is either
Riemannian or Lorentzian. Geometric quantities in $N$ will be denoted by
$(\bar g_{\alpha\beta}),(\riema \alpha\beta\gamma\delta)$, etc., and those in $M$ by $(g_{ij}), (\riem
ijkl)$, etc. Greek indices range from $0$ to $n$ and Latin from $1$ to $n$; the
summation convention is always used. Generic coordinate systems in $N$ resp.
$M$ will be denoted by $(x^\alpha)$ resp. $(\x^i)$. Covariant differentiation will
simply be indicated by indices, only in case of possible ambiguity they will be
preceded by a semicolon, i.e., for a function $u$ in $N$, $(u_\alpha)$ will be the
gradient and
$(u_{\alpha\beta})$ the Hessian, but e.g., the covariant derivative of the curvature
tensor will be abbreviated by $\riema \alpha\beta\gamma{\delta;\e}$. We also point out that
\begin{equation}
\riema \alpha\beta\gamma{\delta;i}=\riema \alpha\beta\gamma{\delta;\e}x_i^\e
\end{equation}
with obvious generalizations to other quantities.

Let $M$ be a \tit{spacelike} hypersurface, i.e. the induced metric is Riemannian,
with a differentiable normal $\n$. We define the signature of $\n$, $\sigma=\s(\n)$, by
\begin{equation}
\sigma=\bar g_{\alpha\beta}\n^\alpha\n^\beta=\spd \n\n.
\end{equation}
In case $N$ is Lorentzian, $\sigma=-1$, and $\n$ is time-like.

In local coordinates, $(x^\alpha)$ and $(\x^i)$, the geometric quantities of the
spacelike hypersurface $M$ are connected through the following equations
\begin{equation}\lae{1.3}
x_{ij}^\alpha=-\sigma h_{ij}\n^\alpha
\end{equation}
the so-called \tit{Gau{\ss} formula}. Here, and also in the sequel, a covariant
derivative is always a \tit{full} tensor, i.e.

\begin{equation}
x_{ij}^\alpha=x_{,ij}^\alpha-\ch ijk x_k^\alpha+\cha \beta\gamma\alpha x_i^\beta x_j^\gamma.
\end{equation}
The comma indicates ordinary partial derivatives.

In this implicit definition the \tit{second fundamental form} $(h_{ij})$ is taken
with respect to $-\sigma\n$.

The second equation is the \tit{Weingarten equation}
\begin{equation}
\n_i^\alpha=h_i^k x_k^\alpha,
\end{equation}
where we remember that $\n_i^\alpha$ is a full tensor.

Finally, we have the \tit{Codazzi equation}
\begin{equation}
h_{ij;k}-h_{ik;j}=\riema\alpha\beta\gamma\delta\n^\alpha x_i^\beta x_j^\gamma x_k^\delta
\end{equation}
and the \tit{Gau{\ss} equation}
\begin{equation}
\riem ijkl=\sigma \{h_{ik}h_{jl}-h_{il}h_{jk}\} + \riema \alpha\beta\gamma\delta x_i^\alpha x_j^\beta x_k^\gamma
x_l^\delta.
\end{equation}
Here, the signature of $\n$ comes into play.

Now, let us assume that $N$ is a globally hyperbolic Lorentzian manifold with a
\tit{compact} Cauchy surface. Then
$N$ is a topological product $\R[]\times \mc S_0$, where $\mc S_0$ is a
compact Riemannian manifold, and there exists a Gaussian coordinate system
$(x^\alpha)$, such that the metric in $N$ has the form 
\begin{equation}\lae{1.8b}
d\bar s_N^2=e^{2\psi}\{-{dx^0}^2+\sigma_{ij}(x^0,x)dx^idx^j\},
\end{equation}
where $\sigma_{ij}$ is a Riemannian metric, $\psi$ a function on $N$, and $x$ an
abbreviation for the spacelike components $(x^i)$,

We also assume that
the coordinate system is \tit{future oriented}, i.e. the time coordinate $x^0$
increases on future directed curves. Hence, the \tit{contravariant} time-like
vector $(\x^\alpha)=(1,0,\dotsc,0)$ is future directed as is its \tit{covariant}
version
$(\x_\alpha)=e^{2\psi}(-1,0,\dotsc,0)$.

Let $M=\graph \fv u\so$ be a spacelike hypersurface
\begin{equation}
M=\set{(x^0,x)}{x^0=u(x),\,x\in\mc S_0},
\end{equation}
then the induced metric has the form
\begin{equation}
g_{ij}=e^{2\psi}\{-u_iu_j+\sigma_{ij}\}
\end{equation}
where $\sigma_{ij}$ is evaluated at $(u,x)$, and its inverse $(g^{ij})=(g_{ij})^{-1}$ can
be expressed as
\begin{equation}\lae{1.10}
g^{ij}=e^{-2\psi}\{\sigma^{ij}+\frac{u^i}{v}\frac{u^j}{v}\},
\end{equation}
where $(\sigma^{ij})=(\sigma_{ij})^{-1}$ and
\begin{equation}\lae{1.11}
\begin{aligned}
u^i&=\sigma^{ij}u_j\\
v^2&=1-\sigma^{ij}u_iu_j\equiv 1-\abs{Du}^2.
\end{aligned}
\end{equation}
Hence, $\graph u$ is spacelike if and only if $\abs{Du}<1$.

The covariant form of a normal vector of a graph looks like
\begin{equation}
(\n_\alpha)=\pm v^{-1}e^{\psi}(1, -u_i).
\end{equation}
and the contravariant version is
\begin{equation}
(\n^\alpha)=\mp v^{-1}e^{-\psi}(1, u^i).
\end{equation}
Thus, we have
\br Let $M$ be spacelike graph in a future oriented coordinate system. Then, the
contravariant future directed normal vector has the form
\begin{equation}
(\n^\alpha)=v^{-1}e^{-\psi}(1, u^i)
\end{equation}
and the past directed
\begin{equation}\lae{1.15}
(\n^\alpha)=-v^{-1}e^{-\psi}(1, u^i).
\end{equation}
\er

In the Gau{\ss} formula \re{1.3} we are free to choose the future or past directed
normal, but we stipulate that in general we use the past directed normal unless otherwise stated.

Look at the component $\alpha=0$ in \re{1.3}, then we obtain in view of \re{1.15}
\begin{equation}\lae{1.16}
e^{-\psi}v^{-1}h_{ij}=-u_{ij}-\cha 000\mspace{1mu}u_iu_j-\cha 0i0
\mspace{1mu}u_j-\cha 0j0\mspace{1mu}u_i-\cha ij0.
\end{equation}
Here, the covariant derivatives a taken with respect to the induced metric of
$M$, and
\begin{equation}\lae{1.18}
-\cha ij0=e^{-\psi}\bar h_{ij},
\end{equation}
where $(\bar h_{ij})$ is the second fundamental form of the hypersurfaces
$\{x^0=\const\}$.

\section{The Beltrami map}\las 2

Let $\R[n+1,1]$ be the $(n+2)$-dimensional Minkowski space with points $x=(x^a)$, $0\le a\le n+1$, where $x^0$ is the time function.

The submanifolds
\begin{equation}
\Hh=\set{x\in \R[n+1,1]}{\spd xx=-1, \,x^0>0}
\end{equation}
and
\begin{equation}
N=\set{x\in \R[n+1,1]}{\spd xx=1}
\end{equation}
are spaces of constant curvature. $\Hh$ is the $(n+1)$-dimensional hyperbolic space with constant curvature $K=-1$, and $N$ is a Lorentzian manifold with constant curvature $K_N=1$, the de Sitter space-time.

$N$ is globally hyperbolic, as can be seen by introducing polar coordinates in the Euclidean part of the Minkowski space such that the metric in $\R[n+1,1]$ is expressed as
\begin{equation}\lae{2.3}
d\bar s^2=-{dx^0}^2+dr^2+r^2\s_{ij}d\xi^i d\xi^j,
\end{equation}
where $\s_{ij}$ is the metric in $S^n$.

Then $N$ is the embedding
\begin{equation}\lae{2.4}
N=\set{(x^0,r,\xi^i)}{r=\sqrt{1+\abs{x^0}^2},\, x^0\in \R[],\,\xi\in S^n},
\end{equation}
i.e., $N=\R[]\times S^n$ topologically and
\begin{equation}\lae{2.5}
ds_N^2=e^{2\psi}\{-d\tau^2+\s_{ij}d\xi^i d\xi^j\},
\end{equation}
where
\begin{equation}\lae{2.6}
\tau=\int_0^{x^0}\frac1{1+t^2}\q\wed\q\psi=\tfrac12\log (1+\abs{x^0}^2).
\end{equation}

Notice that $N$ is simply connected, since $n\ge 2$.

Let us analyze a special representation of $\Hh$ over the unit ball $B_1(0)\su \R$ in some detail.

\bl\lal{2.1}
Let $\pi$ be the so-called Beltrami map
\begin{equation}
\begin{aligned}
\pi:\Hh&\ra B_1(0)\su\R\\
(x^0,x^i)&\ra y=\big(\tfrac{x^i}{x^0}\big).
\end{aligned}
\end{equation}
Then $\pi$ is a diffeomorphism such that, after introducing Euclidean polar coordinates $(r,\xi)$ in $B_1(0)$, the hyperbolic metric can be expressed as
\begin{equation}\lae{2.8}
d\bar s^2=\frac{r^2}{1-r^2}\{\frac1{r^2(1-r^2)}dr^2+\s_{ij}d\xi^i d\xi^j\}
\end{equation}
or, if we define $\tau$ by
\begin{equation}
d\tau=\frac1{r\sqrt{1-r^2}}dr,
\end{equation}
\begin{equation}\lae{2.10}
\begin{aligned}
d\bar s^2&=\frac{r^2}{1-r^2}\{d\tau^2+\s_{ij}d\xi^i d\xi^j\}\\
&\equiv e^{2\psi}\{d\tau^2+\s_{ij}d\xi^i d\xi^j\};
\end{aligned}
\end{equation}
$\tau$ is uniquely determined up to an integration constant.
\el

\bp
Writing the points in $\Hh$ in the form $(x^0,z)$ such that
\begin{equation}
-\abs{x^0}^2+\abs z^2=-1
\end{equation}
we deduce
\begin{equation}
x^0=\frac1{\sqrt{1-\abs y^2}},\q y=\frac z{x^0},
\end{equation}
hence
\begin{equation}
\pi^{-1}(y)=\big(\frac1{\sqrt{1-\abs y^2}}, \frac y{\sqrt{1-\abs y^2}}\big)
\end{equation}
is a bijective mapping from $B_1(0)$ onto $\Hh$.

In polar coordinates $(y^\al)=(r,\xi^i)$, $1\le \al\le n+1$,
\begin{equation}\lae{2.14}
x=\pi^{-1}(y)=\big(\frac1{\sqrt{1-r^2}},\frac{r\xi}{\sqrt{1-r^2}}\big),\q\abs\xi=1,
\end{equation}
and the form \re{2.8} of the hyperbolic metric can be deduced from
\begin{equation}
\bar g_{\al\bet}=\spd{x_\al}{x_\bet},\q 1\le \al,\bet\le n+1.\qedhere
\end{equation}
\ep

Now, let us denote the coordinates $(\tau,\xi^i)$ as usual by $(x^\al)$, $0\le \al\le n$, $\tau=x^0$, and let $(\bar g_{\al\bet})$ be the metric in \re{2.10}.

Let $(\tilde g_{\al\bet})$ be the Euclidean metric in $B_1(0)$ in the coordinate system $(\tilde x^\al)=(\tilde\tau,x^i)$, such that
\begin{equation}
\begin{aligned}
d\tilde s^2=\tilde g_{\al\bet}d\tilde x^\al d\tilde x^\bet&=r^2\{d\tilde\tau^2+\s_{ij}d\xi^id\xi^j\}\\
&\equiv e^{2\tilde\psi}\{d\tilde\tau^2+\s_{ij}d\xi^id\xi^j\},
\end{aligned}
\end{equation}
where
\begin{equation}
d\tilde\tau=r^{-1}dr.
\end{equation}

Writing $\tilde\tau =\f(\tau)$ we deduce
\begin{equation}\lae{2.18}
\frac{d\tilde\tau}{d\tau}=\frac{d\tilde\tau}{dr}\frac{dr}{d\tau}=\dot\f=\sqrt{1-r^2}.
\end{equation}

Let $M\su \Hh$ be an arbitrary closed, connected, strictly convex embedded hypersurface, then $M$ is the boundary of a convex body $\hat M$. Without loss of generality we may assume that $x_0=(1,0)=\pi^{-1}(0)$ is an interior point of $\hat M$. Then $M$ can be written as a graph in geodesic polar coordinates centered at $x_0$, \cf, e.g., \cite[Section 4]{cg96}, or equivalently, as a graph over $S^n$ in the coordinates $(\tau,x^i)$
\begin{equation}
M=\graph u=\set{\tau=u(x)}{x\in S^n}.
\end{equation}

Because of \re{2.18} $M$ can also be viewed as a graph $\tilde M$ in $B_1(0)$ with respect to the Euclidean metric
\begin{equation}
\tilde M=\graph\tilde u=\set{\tilde\tau=\tilde u(x)}{x\in S^n}.
\end{equation}

Then we derive
\begin{equation}
\tilde u=\f(u).
\end{equation}

Denote by $(g_{ij},h_{ij})$ \resp $(\hat g_{ij},\hat h_{ij})$ the metric and second fundamental form of $M$ with respect to the ambient metrics $(\bar g_{\al\bet})$ \resp $(e^{-2\psi}\bar g_{\al\bet})$, and similarly let $(\tilde g_{ij},\tilde h_{ij})$ \resp $(\hat{\tilde g}_{ij},\hat{\tilde h}_{ij})$ be the geometric quantities of $\tilde M$ with respect to the ambient metrics $(\tilde g_{\al\bet})$ \resp $(e^{-2\tilde\psi}\tilde g_{\al\bet})$.

Then we have
\begin{equation}
h_{ij}e^{-\psi}=\hat h_{ij}+\frac{d\psi}{d\tau}v^{-1}\hat g_{ij}
\end{equation}
and
\begin{equation}\lae{2.23}
\tilde h_{ij}e^{-\tilde\psi}=\hat {\tilde h}_{ij}+\frac{d\tilde\psi}{d\tilde\tau}\tilde v^{-1}\hat {\tilde g}_{ij},
\end{equation}
where
\begin{equation}
v^2=1+\s^{ij}u_iu_j
\end{equation}
\begin{equation}
\tilde v^2=1+\s^{ij}\tilde u_i\tilde u_j,
\end{equation}
see e.g., \cite[Prop. 12.2.11]{cg:aii}.

Moreover, there holds
\begin{equation}
\hat g_{ij}=u_iu_j+\s_{ij}
\end{equation}
and
\begin{equation}
\hat h_{ij}=-u_{ij}v^{-1},
\end{equation}
where $(u_{ij})$ is the Hessian of $u$ with respect to the metric $(\s_{ij})$. Analogue formulas are valid for $\hat{\tilde g}_{ij}$ and $\hat{\tilde h}_{ij}$.  

Hence we deduce
\begin{equation}
h_{ij}e^{-\psi}=-u_{ij}v^{-1}+\frac{d\psi}{d\tau}v^{-1}\hat g_{ij}
\end{equation}
and
\begin{equation}
\tilde h_{ij}e^{-\tilde\psi}=-\tilde u_{ij}\tilde v^{-1}+\frac{d\tilde\psi}{d\tilde\tau} \tilde v^{-1}\hat{\tilde g}_{ij}.
\end{equation}

Using the relations
\begin{equation}
\tilde u_i=\dot\f u_i,
\end{equation}
and
\begin{equation}
\tilde u_{ij}=\dot \f u_{ij}+\Ddot\f u_i u_j,
\end{equation}
where
\begin{equation}
\dot\f=\sqrt{1-r^2},\q \Ddot\f =-r^2,
\end{equation}
we obtain after some elementary calculations
\begin{equation}
\tilde h_{ij}\tilde v=(1-r^2)h_{ij}v,
\end{equation}
i.e., $\tilde M$ is also strictly convex.

Moreover, let $u^i=\s^{ij}u_j$, then
\begin{equation}
\begin{aligned}
\tilde g^{ij}&=r^{-2}\{\s^{ij}-\frac{\dot\f^2}{\tilde v^2}u^iu^j\}\\
&\ge r^{-2}\{\s^{ij}-\frac{1}{v^2}u^iu^j\},
\end{aligned}
\end{equation}
since
\begin{equation}
\dot\f^2=1-r^2<1
\end{equation}
and we conclude
\begin{equation}
\tilde h^j_i\tilde v\ge h^j_iv,
\end{equation}
hence,
\begin{equation}
\tilde h^j_i\ge h^j_i\frac v{\tilde v}\ge h^j_i.
\end{equation}

Note also, that in points where $Du=0$ there holds
\begin{equation}
\tilde h^j_i=h^j_i,
\end{equation}
i.e., the principal curvatures are then identical.

Thus, we have proved
\bl\lal{2.2}
Let $M\su\Hh$ be a closed, connected, strictly convex hypersurface, then the Beltrami map $\pi$ maps $M$ onto a closed strictly convex hypersurface $\tilde M\su B_1(0)$. Moreover, expressing the normal vectors $\nu$ \resp $\tilde \nu$ of $M$ \resp $\tilde M$ in the common coordinate system $(\tau, \xi^i)$ yields that they are collinear.
\el

\bp
Only the last statement needs a verification. Up to a positive factor the covariant normal vector $(\nu_\al)$ has the form
\begin{equation}\lae{2.39}
(\nu_\al)=(1,-u_i)
\end{equation}
and $(\tilde \nu_\al)$, in the coordinate system $(\tau,\xi^i)$,
\begin{equation}\lae{2.40}
\begin{aligned}
(\tilde\nu_\al)=(\frac{d\tilde\tau}{d\tau},-\tilde u_i)=(\dot\f,-\dot\f u_i)=\dot\f(1,-u_i).\qedhere
\end{aligned}
\end{equation}
\ep

\br\lar{2.3}
The results of the preceding lemma can also be applied to a local embedding of a strictly convex hypersurface $M$ that can be represented as a graph in geodesic polar coordinates centered in the Beltrami point $(1,0)$ regardless which side of $M$ the Beltrami point is facing.
\er

\section{Hadamard's theorem in hyperbolic space}\las{3}

\bt\lat{3.1}
Let $M_0$ be a compact, connected $n$-dimensional manifold and
\begin{equation}
x:M_0\ra \Hh
\end{equation}
a strictly convex immersion of class $C^2$, i.e., the second fundamental form with respect to any normal is always (locally) invertible, then the immersion is actually an embedding and $M=x(M_0)$ a strictly convex hypersurface that bounds a strictly convex body $\hat M\su\Hh$. $M$ and $M_0$ are moreover diffeomorphic to $S^n$ and orientable.
\et

\bp
Since we shall again employ the Beltrami map, we consider $\Hh$ as a hypersurface in $\R[n+1,1]$ and $M$ as a codimension $2$ immersed submanifold in $\R[n+1,1]$, i.e.,
\begin{equation}
x:M_0\ra \R[n+1,1].
\end{equation}

The Gaussian formula for $M$ then looks like
\begin{equation}\lae{3.3}
x_{ij}=g_{ij}x-h_{ij}\tilde x,
\end{equation}
where $g_{ij}$ is the induced metric, $h_{ij}$ the second fundamental form of $M$ considered as a hypersurface in $\Hh$, and $\tilde x$ is the representation of the (exterior\footnote{Notice that for any closed, connected immersed hypersurface in $\Hh$ an exterior normal vector can be unambiguously defined.}) normal vector $\nu=(\nu^\al)$ of $M$ in $T(\Hh)$ as a vector in $T(\R[n+1,1])$.

Without loss of generality we may assume that the \tit{Beltrami point} $(1,0)$ does not belong to $M$, since the isometries of $\Hh$ act transitively. Let $\pi$ be the Beltrami map such that
\begin{equation}
\pi(1,0)=0\in\R
\end{equation}
and denote by $\f$ its inverse
\begin{equation}\lae{3.6}
\f=\pi^{-1}:\R\ra \Hh\su \R[n+1,1].
\end{equation}

Corresponding to the immersion $x=x(\xi)$ we then have an immersion $y=\pi\circ x$
\begin{equation}
y:M_0\ra \R.
\end{equation}

Let $\tilde M\su\R$ be its image and $\tilde g_{ij}$ $\tilde h_{ij}$, $\tilde \nu=(\tilde\nu^\al)$ its geometric quantities. We shall prove that $\tilde M$ is an immersed, closed strictly convex hypersurface and hence an embedded hypersurface, due to Hadamard's theorem, \cf \cite{sacksteder}.

In view of the relation
\begin{equation}\lae{3.7}
x=\f\circ y
\end{equation}
we shall then deduce that $x=x(\xi)$ is an embedding.

The inverse Beltrami map $\f$ provides an embedding of $\Hh$ in $\R[n+1,1]$, i.e., we have the Gaussian formula
\begin{equation}\lae{3.8}
\f_{\al\bet}=\bar g_{\al\bet}\f,
\end{equation}
where $\bar g_{\al\bet}$ is the induced metric
\begin{equation}
\bar g_{\al\bet}=\spd{\f_\al}{\f_\bet}.
\end{equation}

Differentiating \re{3.7} covariantly with respect to the metric $g_{ij}$ of $M$ we obtain
\begin{equation}\lae{3.10}
\begin{aligned}
x_{ij}&=\f_\al y^\al_{ij}+\f_{\al\bet}y^\al_iy^\bet_j\\
&=\f_\al y^\al_{ij}+ \bar g_{\al\bet}y^\al_iy^\bet_j\f,
\end{aligned}
\end{equation}
in view of \re{3.8}. 

Indicate covariant derivatives with respect to the metric $\tilde g_{ij}$ by a preceding semicolon such that
\begin{equation}
y_{;ij}=-\tilde h_{ij}\tilde \nu,
\end{equation}
then
\begin{equation}
y_{ij}=y_{;ij}-\{\ch ijk-\tilde\C^k_{ij}\}y_k,
\end{equation}
hence, we derive from \re{3.10}
\begin{equation}\lae{3.13}
x_{ij}=\f_\al y^\al_{;ij}-\{\ch ijk-\tilde\C^k_{ij}\}y^\al_k\f_\al+\bar g_{\al\bet}y^\al_iy^\bet_j\f.
\end{equation}

Let $\nu=(\nu^\al)$ be the exterior normal of $M$ expressed in the coordinates $(y^\al)$ of $\R$, then the representation $\tilde x$ of $\nu$ in $T(\R[n+1,1])$is given by
\begin{equation}
\tilde x=\f_\al\nu^\al
\end{equation}
as can be easily checked.

From \re{3.3} and \re{3.13} we then deduce
\begin{equation}
h_{ij}=-\spd{x_{ij}}{\f_\ga\nu^\ga}=\tilde h_{ij}\bar g_{\al\bet}\nu^\al\tilde\nu^\bet,
\end{equation}
where we used
\begin{equation}
x_k=\f_\al y^\al_k
\end{equation}
and
\begin{equation}\lae{3.17}
\spd{x_k}{\tilde x}=0.
\end{equation}

Thus, it remains to prove that
\begin{equation}
\bar g_{\al\bet}\nu^\al\tilde\nu^\bet\ne 0\qq\tup{in}\; M_0.
\end{equation}

So far we haven't used the fact $0\notin \tilde M$, or equivalently, $(1,0)\notin M$, but now we introduce polar coordinates $(y^\al)$ in $\R$ such that $y^0=r$ and distinguish two cases
\begin{equation}\lae{3.19}
\spd{\frac{\pa}{\pa r}}\nu=0
\end{equation}
and
\begin{equation}\lae{3.20}
\spd{\frac{\pa}{\pa r}}\nu\ne 0,
\end{equation}
where the metric is the one in $\Hh$.

In Euclidean polar coordinates $(y^\al)=(r,\xi^i)$ the hyperbolic metric $(\bar g_{\al\bet})$ has been expressed in \re{2.8}. Hence, if \re{3.19} is valid, we deduce
\begin{equation}
\nu^0=0,
\end{equation}
and infer further, in view of \re{3.17},
\begin{equation}
\begin{aligned}
0&=\spd{x_i}{\tilde x}=\spd{\f_\al}{\f_\bet}y^\al_i\nu^\bet\\
&=\bar g_{\al\bet}y^\al_i\nu^\bet=\frac{r^2}{1-r^2}\s_{kj}y^k_i\nu^j,
\end{aligned}
\end{equation}
from which we conclude that
\begin{equation}
\tilde g_{\al\bet}y^\al_i\nu^\bet=0,
\end{equation}
where $\tilde g_{\al\bet}$ is the Euclidean metric expressed in polar coordinates.

Thus, $\nu$ and $\tilde \nu$ are collinear, if \re{3.19} is valid.

On the other hand, if the assumption \re{3.20} is satisfied, then $M$, or more precisely, a local embedding of $M_0$ can be written as a graph in polar coordinates, i.e., we are in the situation where the results of \rl{2.2} and the equations \re{2.39}, \re{2.40} can be applied locally, \cf \rr{2.3}, and we deduce again that $\nu$, $\tilde\nu$ are collinear.

Therefore, Hadamard's theorem yields that the immersion is actually an embedding and that $M_0$, and hence $M$, is diffeomorphic to $S^n$.
\ep

\section{The Gau{\ss} maps}\las{4}

Let $M\su\Hh$ be a closed, connected, strictly convex hypersurface given by an embedding
\begin{equation}
x:M_0\ra M.
\end{equation}

Considering $M$ as a codimension $2$ submanifold of $\R[n+1,1]$ such that
\begin{equation}
x_{ij}=g_{ij}x-h_{ij}\tilde x,
\end{equation}
where $\tilde x\in T_x(\R[n+1,1])$ represents the exterior normal vector $\nu\in T_x(\Hh)$, we want to prove that the mapping
\begin{equation}\lae{4.3}
\tilde x:M_0\ra N
\end{equation}
is an embedding of a strictly convex, closed, spacelike hypersurface $\tilde M$. We call this mapping the \tit{Gau{\ss} map} of $M$.

First, we shall show that the Gau{\ss} map is injective. To prove this result we need the following lemma.
\bl\lal{4.1}
Let $M\su\Hh$ be a closed, connected, strictly convex hypersurface and denote by $\hat M$ its (closed) convex body. Let $x\in M$ be fixed and $\tilde x$ be the corresponding outward normal vector, then
\begin{equation}\lae{4.4}
\spd y{\tilde x}\le 0\qq\A\, y\in \hat M
\end{equation}
and also strictly less than $0$ unless $y=x$.

The preceding inequality also characterizes the points in $\hat M$, namely, let $y\in \Hh$ be such that
\begin{equation}\lae{4.5b}
\spd y{\tilde x}\le 0\qq\A\,x\in M,
\end{equation}
then $y\in \hat M$.
\el

\bp
\cq{\re{4.4}}\q   Let $y\in \tup{int}\, \hat M$ be arbitrary and let $z=z(t)$, $0\le t\le d$ be the unique geodesic in $\Hh$ connecting $y$ and $x$ such that
\begin{equation}
z(0)=x\q\wed\q z(d)=y
\end{equation}
parametrized by arc length.

Viewing $z$ as a curve in $\R[n+1,1]$ the geodesic equation has the form
\begin{equation}
\Ddot z\equiv \tfrac{D}{dt}\dot z =z.
\end{equation}
If the coordinate system in $\R[n+1,1]$ is Euclidean, the covariant derivatives are just ordinary derivatives.

It is well-known that the geodesic $z$ is contained in $\hat M$ and that
\begin{equation}\lae{4.8}
\spd{\dot z(0)}{\tilde x}<0;
\end{equation}
notice that, after introducing geodesic polar coordinates in $\Hh$ centered in $y$, we have
\begin{equation}
\spd{\dot z(0)}{\tilde x}=-\spd{\frac{\pa}{\pa r}}\nu
\end{equation}
and hence is strictly negative, \cf \cite[Section 4]{cg96}.

Thus, $\f(t)=\spd{z(t)}{\tilde x}$ satisfies the initial value problem
\begin{equation}\lae{4.9}
\Ddot \f=\f,\q \f(0)=0,\q \dot\f(0)<0,
\end{equation}
and is therefore equal to
\begin{equation}
\f(t)=-\lam \sinh t,\q \lam>0,
\end{equation}
i.e.,
\begin{equation}
\f(t)<0\qq\A\,t>0.
\end{equation}

Now, let $y\in M$, $y\ne x$, be arbitrary, and consider a sequence $z_k$ of geodesics parametrized in the interval $0\le t\le 1$, such that
\begin{equation}
z(0)=x\q\wed\q z_k(1)\ra y,
\end{equation}
where $z_k(1)\in \tup{int}\, \hat M$.

The geodesics $z_k$ converge to a geodesic $z$ connecting $x$ and $y$. If
\begin{equation}
\spd{\dot z(0)}{\tilde x}<0,
\end{equation}
then the previous arguments are valid yielding
\begin{equation}
\spd y{\tilde x}<0.
\end{equation}

On the other hand, the alternative
\begin{equation}
\spd y{\tilde x}=0
\end{equation}
leads to a contradiction, since then the geodesic $z$ would be part of the tangent space $T_x(M)$ which is impossible, \cf the considerations in  \cite{cg96} after the equation $(4.17)$.

\cvm
\cq{$y\in \hat M$}\q Suppose now that $y\in \Hh$ satisfies \re{4.5b}, and assume by contradiction that $y\in \cpl \hat M$. Pick an arbitrary $x_0\in \tup{int}\msp \hat M$ and let $z=z(t)$, $0\le t\le d$,  be the geodesic joining $x_0$ and $y$ parameterized by arc length, such that $z(0)=x_0$ and $z(d)=y$. The geodesic intersects $M$ in a unique point $x$, $x=z(t_1)$, $0<t_1<d$.

Define
\begin{equation}
\f(t)=\spd {z(t)}{\tilde x}
\end{equation}
and let $0\le t_0\le d$ be such that
\begin{equation}
\f(t_0)=\sup\set{\f(t)}{0\le t\le d}.
\end{equation}

We now distinguish two cases. First, we assume $\f(t_0)>0$, then there must hold $0<t_0<d$ and $\dot \f(t_0)=0$. Thus $\f$ satisfies the initial value problem
\begin{equation}
\Ddot \f=\f,\q \f(t_0)>0,\q \dot\f(t_0)=0,
\end{equation}
and must therefore be equal to
\begin{equation}
\f(t)=\lam \cosh(t-t_0), \q \lam >0,
\end{equation}
which is a contradiction, since $\f(0)<0$.

Hence, we must have $\f(t_0)=0$ and we may choose $t_0=t_1$, i.e., there holds $\dot\f(t_1)=0$, which is a contradiction too, because of  the inequality \re{4.8}, which now reads $\dot\f(t_1)>0$.

Therefore we have proved $y\in\hat M$.
\ep

\bt\lat{4.2}
Let $x:M_0\ra M\su \Hh$ be the embedding of a closed, connected, strictly convex hypersurface, then the Gau{\ss} map defined in \re {4.3} is injective, where we identify $\R[n+1,1]$ with its individual tangent spaces.
\et

\bp
We again assume $M$ to be a codimension $2$ submanifold in $\R[n+1,1]$. Suppose there would be two points $p_1\ne p_2$ in $M_0$ such that
\begin{equation}
\tilde x(p_1)=\tilde x(p_2),
\end{equation}
then the function
\begin{equation}
\f(y)=\spd y{\tilde x(p_1)}
\end{equation}
would vanish in the points $x(p_1)$ as well as $x(p_2)$ contrary to the results of \rl{4.1}.
\ep

\bl
As a submanifold of codimension $2$ $M$ satisfies the Wein\-gar\-ten equations
\begin{equation}\lae{4.18}
\tilde x_i=g^k_ix_k
\end{equation}
for the normal $\tilde x$ and also
\begin{equation}
x_i=g^k_ix_k
\end{equation}
for the normal $x$.
\el

\bp
We only have to prove the non-trivial Weingarten equation.

First we infer from
\begin{equation}
\spd x{\tilde x}=0
\end{equation}
that
\begin{equation}
0=\spd{x_i}{\tilde x}+\spd x{\tilde x_i}=\spd x{\tilde x_i}.
\end{equation}
Furthermore, there holds
\begin{equation}\lae{4.22}
0=\spd{\tilde x}{\tilde x_i},
\end{equation}
since $\spd {\tilde x}{\tilde x}=1$. Hence, we deduce
\begin{equation}\lae{4.23}
\tilde x_i=a^k_ix_k.
\end{equation}

Differentiating the relation $\spd{x_j}{\tilde x}=0$ covariantly we obtain
\begin{equation}\lae{4.24}
\spd {\tilde x_j}{x_i}=h_{ij}
\end{equation}
and we infer \re{4.18} in view of \re{4.23}.
\ep

We can now prove
\bt\lat{4.4}
Let $x:M_0\ra M\su \Hh$ be a closed, connected, strictly convex hypersurface of class $C^m$, $m\ge 3$, then the Gau{\ss} map $\tilde x$ in \re{4.3} is the embedding of a closed, spacelike, achronal, strictly convex hypersurface $\tilde M\su N$ of class $C^{m-1}$. 

Viewing $\tilde M$ as a codimension $2$ submanifold in $\R[n+1,1]$, its Gaussian formula is
\begin{equation}\lae{4.25}
\tilde x_{ij}=-\tilde g_{ij}\tilde x+\tilde h_{ij} x,
\end{equation}
where $\tilde g_{ij}$, $\tilde h_{ij}$ are the metric and second fundamental form of the hypersurface $\tilde M\su N$, and $x=x(\xi)$ is the embedding of $M$ which also represents the future directed normal vector of $\tilde M$. The second fundamental form $\tilde h_{ij}$ is defined with respect to the future directed normal vector, where the time orientation of $N$ is inherited from $\R[n+1,1]$.

The second fundamental forms of $M$, $\tilde M$ and the corresponding principal curvatures $\ka_i$, $\tilde \ka_i$ satisfy
\begin{equation}\lae{4.26}
h_{ij}=\tilde h_{ij}=\spd{\tilde x_i}{x_j}
\end{equation}
and
\begin{equation}\lae{4.27}
\tilde \ka_i=\ka_i^{-1}.
\end{equation}
\et

\bp
(i) From the Weingarten equation \re{4.18} we infer
\begin{equation}\lae{4.28}
\tilde g_{ij}=\spd{\tilde x_i}{\tilde x_j}=h^k_ih_{kj}
\end{equation}
is positive definite, hence $\tilde x=\tilde x(\xi)$ is an embedding of a closed, connected spacelike hypersurface, where we also used \rt{4.2}.

Since $N$ is simply connected, we conclude further that $\tilde M$ is achronal and thus can be written as a graph over the Cauchy hypersurface $\{0\}\ti S^n$ which we identify with $S^n$
\begin{equation}
\tilde M=\graph \fv{\tilde u}{S^n}\su N,
\end{equation}
\cf \cite[p. 427]{bn} and \cite[Prop. 2.5]{cg:indiana}.

\cvm
(ii) The pair $(x,\tilde x)$ satisfies
\begin{equation}\lae{4.30}
\spd x{\tilde x}=0
\end{equation}
and we claim that $x$ is the future directed normal vector of $\tilde M$ in $\tilde x$, where as usual we identify the normal vector $\tilde \nu=(\tilde\nu^\al)\in T_{\tilde x}(N)$ with its embedding in $T_{\tilde x}(\R[n+1,1])$.

Differentiating \re{4.30} covariantly and using the fact that $\tilde x$ is a normal vector for $M$ we deduce
\begin{equation}\lae{4.31}
0=\spd x{\tilde x_i},
\end{equation}
i.e., $\tilde x$ and $x$ span the normal space of the codimension $2$ submanifold $\tilde M$. By the very definition of $\Hh$ $x$ is a future directed vector in $\R[n+1,1]$.

Let us define the second fundamental form $\tilde h_{ij}$ of $\tilde M\su N$ with respect to the future directed normal vector $\tilde \nu \in T_{\tilde x}(N)$, then the codimension $2$ Gaussian formula is exactly \re{4.25} because of \re{4.31}.  

Differentiating the Weingarten equation \re{4.18} covariantly with respect to the metric $\tilde g_{ij}$ and indicating the covariant derivatives with respect to $\tilde g_{ij}$ by a semi-colon and those with respect to $g_{ij}$ simply by indices, we obtain
\begin{equation}
\tilde x_{;ij}=h^k_{i;j}x_k+h^k_ix_{;kj}
\end{equation}
and we deduce further
\begin{equation}\lae{4.33}
\tilde h_{ij}=-\spd{\tilde x_{;ij}}x=-h^k_i\spd{x_{kj}}x=h^k_ig_{kj}=h_{ij}.
\end{equation}

On the other hand, we infer from \re{4.31}
\begin{equation}
\tilde h_{ij}=-\spd{\tilde x_{;ij}}x=\spd{\tilde x_i}{x_j}
\end{equation}
which proves \re{4.26}.

The last relation \re{4.27} follows from \re{4.33} and \re{4.28}. 
\ep

We can also define a Gau{\ss} map from strictly convex, connected, spacelike hypersurfaces $\tilde M\su N$ into $\Hh$ such that the two Gau{\ss} maps are inverse to each other.

\bt\lat{4.5}
Let $\tilde M\su N$ be a closed, connected, spacelike, strictly convex, embedded hypersurface of class $C^m$, $m\ge 3$, such that, when viewed as a codimension $2$ submanifold in $\R[n+1,1]$, its Gaussian formula is
\begin{equation}
\tilde x_{ij}=-\tilde g_{ij}\tilde x+\tilde h_{ij}x,
\end{equation}
where $\tilde x=\tilde x(\xi)$ is the embedding, $x$ the future directed normal vector, and $\tilde g_{ij}$, $\tilde h_{ij}$ the induced metric and the second fundamental form of the hypersurface in $N$. Then we define the Gau{\ss} map as $x=x(\xi)$
\begin{equation}
x:\tilde M\ra \Hh\su \R[n+1,1].
\end{equation}
The Gau{\ss} map is the embedding of a closed, connected, strictly convex hypersurface $M$ in $\Hh$.

Let $g_{ij}$, $h_{ij}$ be the induced metric and second fundamental form of $M$, then, when viewed as a codimension $2$ submanifold, $M$ satisfies the relations
\begin{equation}\lae{4.37}
x_{ij}=g_{ij}x-h_{ij}\tilde x,
\end{equation}
\begin{equation}\lae{4.38}
h_{ij}=\tilde h_{ij}=\spd{x_i}{\tilde x_j},
\end{equation}
and
\begin{equation}\lae{4.39}
\ka _i=\tilde \ka_i^{-1},
\end{equation}
where $\ka_i$, $\tilde \ka_i$ are the corresponding principal curvatures.
\et

\bp
The fact that $x=x(\xi)$ is the immersion of a closed, connected,  strictly convex hypersurface $M$ satisfying the relations \re{4.37}, \re{4.38}, and \re{4.39} follows along the lines of the proof of the previous theorem.

Using \rt{3.1} we then deduce that the immersion is an embedding.
\ep

Combining the two theorems, looking especially at the Gaussian formulas \re{4.25} and \re{4.37}, we immediately conclude that the Gau{\ss} maps are inverse to each other, i.e., if we start with a closed, strictly convex hypersurface $M\su\Hh$, apply the Gau{\ss} map to obtain a spacelike,  strictly convex hypersurface $\tilde M\su N$, and then apply the second Gau{\ss} map, then we return to $M$ with a pointwise equality.

Denoting the two Gau{\ss} maps simply by a tilde, this can be expressed in the form
\begin{equation}
x=\tilde{\tilde x},
\end{equation}
or, equivalently, in the form of a commutative diagram
\begin{equation}
\begin{diagram}
\node{M}
      \arrow[2]{e,t}{\tilde{}} \arrow{se,r}{\id}
   \node[2]{\tilde M} \arrow{sw,b}{\tilde{}}\\
\node[2]{M}
\end{diagram}
\end{equation}

Before we give an equivalent characterization of the images of the Gau{\ss} maps, let us show that the images of  strictly convex hypersurfaces by the Gau{\ss} maps are as smooth as the original hypersurfaces.

\bt
Let $M\su \Hh$ be a closed, connected, strictly convex hypersurface of class $C^{m,\al}$, $m\ge 2$, $0\le \al\le 1$, then $\tilde M\su N$, its image under the Gau{\ss} map is also of class $C^{m,\al}$. 

The corresponding regularity result is also valid, if we start with a closed, spacelike, connected, strictly convex hypersurface in $N$ and use the Gau{\ss} map to embed it into $\Hh$.
\et

\bp
We only consider the case when we apply the Gau{\ss} map to $M\su \Hh$. Moreover, without loss of generality we shall also assume that the Beltrami point $(1,0)$ is not part of $M$.

(i) First, let us assume that $m\ge 3$ and $0\le \al\le 1$. The Gau{\ss} map is then of class $C^{m-1,\al}$, i.e., $\tilde M$ is of class $C^{m-1,\al}$. Here, we use the coordinates $(\xi^i)$ for $M$ also as coordinates for $\tilde M$. The metric $\tilde g_{ij}$ and the Christoffel symbols  of $\tilde M$ are then  of class $C^{m-2,\al}$ \resp $C^{m-3,\al}$, while the second fundamental form $\tilde h_{ij}$ is of class $C^{m-2,\al}$, in view of \re{4.26}.

Representing now $\tilde M$ as a graph over $S^n$,
\begin{equation}
\tilde M=\graph\fv u{S^n}
\end{equation}
in conformal coordinates, i.e., we use the coordinates defined in the formulas  \re{2.3} to \re{2.6} denoting them this time, however, by $(\tau, x^i)$ instead of $(\tau,\xi^i)$, since $(\xi^i)$ are supposed to be  given coordinates for $M$. 

Notice that the transformation  $(x^i(\xi^k))$ is a diffeomorphism of class $C^{m-1,\al}$, since the underlying polar coordinates $(x^0,r,x^i)$, defined in \re{2.3}, also cover that part of $\R[n+1,1]$ that contains $M$, due to our assumption at the beginning of the proof. Hence, expressing the Gau{\ss} map in this ambient coordinate system
\begin{equation}
\tilde x(\xi)=(x^0(\xi),r(\xi),x^i(\xi)),
\end{equation}
where 
\begin{equation}
r=\sqrt{1+\abs{x^0}^2},
\end{equation}
we deduce
\begin{equation}
\begin{aligned}
\tilde g_{ij}&=\spd{\tilde x_i}{\tilde x_j}= -x^0_ix^0_j+r_ir_j+r^2\s_{kl}x^k_ix^l_j\\
&=-\frac 1{1+\abs{x^0}^2} x^0_ix^0_j+(1+\abs{x^0}^2)\msp \s_{kl}x^k_ix^k_j
\end{aligned}
\end{equation}
in view of \re{2.3}, proving that the Jacobian $(x^k_i)$ is invertible.

Thus, we conclude that the second fundamental form $\tilde h_{ij}$ expressed in the new coordinates $(x^i)$ is still of class $C^{m-2,\al}$.

We want to express the covariant derivatives $u_{ij}$ of $u$ with respect to the metric $\s_{ij}$ in terms of  $\tilde h_{ij}$ to deduce that $u_{ij}$ is of class $C^{m-1,\al}$, and hence $u\in C^{m,\al}(S^n)$.

To achieve this we define a new metric $\hat g_{\al\bet}$ in the ambient space
\begin{equation}
\hat g_{\al\bet}=e^{-2\psi}\tilde g_{\al\bet},
\end{equation}
where $\tilde g_{\al\bet}$ is the metric in \re{2.5}. 
Let $\hat g_{ij}$, $\hat h_{ij}$ and $\hat \nu$ be the obvious geometric quantities of $\tilde M$ with respect to the new metric, then there holds
\begin{equation}
h_{ij} e^{-\psi}=\hat h_{ij}+\psi_\al\hat\nu^\al \hat g_{ij}
\end{equation}
\cf \re{2.23}, where we already used this formula.

On the other hand, $\hat h_{ij}$ can be expressed in terms of the Hessian $u_{;ij}$ of $u$ with respect to  the metric $\s_{ij}$, namely,
\begin{equation}
\hat h_{ij}=u_{;ij}v^{-1},
\end{equation}
i.e.,
\begin{equation}
h_{ij}e^{-\psi}=u_{;ij} +\psi_\al\hat\nu^\al(-u_iu_j+\s_{ij}),
\end{equation}
hence, $u_{;ij}$ is of class $C^{m-2,\al}$.

\cvm
(ii) The case $m=2$ and $0\le \al\le 1$ follows by approximation and the uniform $C^{2,\al}$-estimates. Notice that the approximating second fundamental forms will converge in $C^0$.
\ep

\bd
(i) Let $M\su \Hh$ be a closed, connected, strictly convex hypersurface, then we define its \tit{polar set} $M^*\su N$ by
\begin{equation}
M^*=\set{y\in N}{\sup_{x\in M}\spd xy=0},
\end{equation}
where the scalar product is the scalar  product in $\R[n+1,1]$ and $x$, $y$ are Euclidean coordinates. 

\cvm
(ii) A similar definition holds, if $M\su N$ is  a spacelike, closed, connected, strictly convex hypersurface $M\su N$, then
\begin{equation}
M^*=\set{y\in \Hh}{\sup_{x\in M}\spd xy=0}.
\end{equation}
\ed

\bt
The polar sets  agree with the images of the Gau{\ss} maps.
\et

\bp
Again we only consider the case $M\su \Hh$.

 In view of \rl{4.1} there holds
 \begin{equation}
\tilde M\su M^*.
\end{equation}

On the other hand, let $y\in M^*$ and $x\in M$ be such that
\begin{equation}\lae{4.53}
\spd xy=0.
\end{equation}
Then we deduce, after introducing local coordinates in $M$,
\begin{equation}\lae{4.54}
\spd{x_i}y=0
\end{equation}
and
\begin{equation}\lae{4.55}
\spd{x_{ij}}y\le 0,
\end{equation}
where the derivatives are covariant derivatives with respect to the induced metric $g_{ij}$ of $M$ being viewed as a codimension $2$ submanifold.

Combining \re{4.53} and \re{4.54} we infer
\begin{equation}
y=\pm \tilde x,
\end{equation}
but because of \re{4.37} and \re{4.55} we deduce $y=\tilde x$.
\ep

Let us conclude the section with
\bt\lat{4.9}
The Gau{\ss} maps provide a bijective relation between the connected, closed, strictly convex hypersurfaces $M\su \Hh$ having the Beltrami point in the interior of their convex bodies and the spacelike, closed,  connected, strictly convex hypersurfaces $\tilde M\su N_+$, where
\begin{equation}
N_+=\set{x\in N}{x^0>0}.
\end{equation}
The geodesic spheres with center in the Beltrami point are mapped onto the coordinate slices $\{x^0=\const\}$.
\et

\bp
(i) Let $M\su \Hh$ be closed, strictly convex such that $p_0\in \tup{int}\msp \hat M$, where $p_0=(1,0,\ldots,0)\in \R[n+1,1]$. According to \rl{2.1}, $\Hh\su \R[n+1,1]$ can be written as the embedding of $\R$ via the inverse $\f=\pi^{-1}$ of the Beltrami map $\pi$, and  $M$ can be represented as $M=\fv{\graph u}{S^n}$ in geodesic polar coordinates centered in $p_0$, or more precisely, in the coordinates $(\tau,\xi^i)$.

A moment's reflection reveals that the Gau{\ss} map of $M$ is given by
\begin{equation}\lae{4.63}
\tilde x=\f_\al\nu^\al,\qq\n \in T_x(\Hh),
\end{equation}
where $\n$ is the exterior normal. In geodesic polar coordinates the normal $\nu$ is given by
\begin{equation}\lae{4.64}
(\nu^\al)=\tilde v e^{-\psi}(1,-u^i),
\end{equation}
where $\tilde v=v^{-1}$ and $u^i=\s^{ij}u_j$; notice that the metric in $\Hh$ is expressed as in \re{2.10}.

Hence, we deduce from \re{2.14}
\begin{equation}
\tilde x^0= \frac{r^2}{(1-r^2)^{3/2}} \tilde v e^{-\psi}=\frac r{\sqrt{1-r^2}}\tilde v>0,
\end{equation}
i.e., $\tilde M\su N_+$.

If $M$ is a geodesic sphere, then we deduce from \re{4.63} and \re{4.64} that it is mapped onto a coordinate slice in $N_+$.

\cvm
(ii) To prove the inverse relation, consider a spacelike, closed,  connected, strictly convex hypersurface $M\su N_+$. Assuming the coordinate system in \re{2.4}, \re{2.5}, $N_+$ can be viewed as the embedding of $\R\sminus \bar B_1(0)$ in $\R[n+1,1]$ via the map 
\begin{equation}\lae{4.66}
x=\f(r,\xi)=(\sqrt{r^2-1},r\xi),\q  \xi\in S^n.
\end{equation}

The Gau{\ss} map from $M$ into $\Hh$ can then be expressed as
\begin{equation}\lae{4.67}
\tilde x=\f_\al\nu^\al,
\end{equation}
where
\begin{equation}
\nu=\tilde ve^{-\psi}(1,u^i)
\end{equation}
is the future directed normal vector in $x\in M$. Let $\tilde M\su \Hh$ be its image. Then we have to show that the Beltrami point $p_0$ is an interior point  of the corresponding convex body, or equivalently, that
\begin{equation}
\spd {p_0}x<0\qq\A\, x\in M,
\end{equation}
in view of the second part of \rl{4.1}. But we immediately deduce
\begin{equation}
\spd{p_0}x=-\sqrt{r^2-1},
\end{equation}
in view of \re{4.66}.

Again we conclude from \re{4.67} that coordinate slices are mapped onto geodesic spheres.
\ep

\section{Curvature flow}\las{5}

Let us now consider the problem of finding a solution of 
\begin{equation}
\fv FM=f(\nu),
\end{equation}
where $F$ is a curvature function defined in the open positive cone $\C_+\su \R[n]$, $0<f$ is a function defined in the normal space of $M$, and $M\su\Hh$ is a closed, connected, strictly convex hypersurface yet to be determined.

Using the results of the previous section, especially \rt{4.4} and \rt{4.5}, we can reformulate the problem equivalently by assuming that $0<f\in C^{2,\al}(N)$, $0<\al<1$, is given and a closed, strictly convex hypersurface $M\su\Hh$ is to be found satisfying
\begin{equation}\lae{5.2}
\fv FM=f(\tilde x)\qq\A\,x\in M,
\end{equation}
where $x\ra \tilde x$ is the Gau{\ss} map corresponding to $M$.

Let $\tilde M\su N$ be the image of $M$ under the Gau{\ss} map, which is identical with the polar $M^*$ of $M$, and let $\tilde F$ be the inverse of $F$, i.e.,
\begin{equation}
\tilde F(\ka_i)=\frac1{F(\ka_i^{-1})},
\end{equation}
then the equation \re{5.2} is equivalent to
\begin{equation}\lae{5.4}
\fv{\tilde F}{\tilde M}=f^{-1}(\tilde x)\qq\A\,\tilde x\in \tilde M,
\end{equation}
in view of \re{4.27}, where now the right-hand side depends on the points $\tilde x\in\tilde M$, and $\tilde M\su N$ is a closed, spacelike, connected strictly convex hypersurface in the de Sitter space $N$ with curvature $K_N=1$.

We solved problems of this kind in \cite[Theorem 0.2]{cg:indiana} assuming barrier conditions and some additional hypotheses. In that paper we denoted the curvature function, the right-hand side and the hypersurface by $F$, $f$, and $M$, which would correspond to the present notation $\tilde F$, $f^{-1}$, $\tilde M$.

Let us stick to the last notation just long enough to formulate the condition for $\tilde F$, namely, we assume that $\tilde F$ is of class $(K^*)$, where the class $(K^*)$ is defined in \cite[Definition 1.6]{cg:indiana} as

\bd
A symmetric curvature function $F\in C^{2.\al}(\C_+)\ii C^0(\bar\C_+)$ of class $(K)$\footnote{For a definition of the class $(K)$ see \cite[Definition 1.1]{cg:indiana}.} is said to be of class $(K^*)$, if there exists $0<\e_0=\e_0(F)$ such that
\begin{equation}\lae{5.5}
\e_0 FH\le F^{ij}h_{ik}h^k_j
\end{equation}
for any positive symmetric tensor $(h_{ij})$. Here $(h_{ij})$ and a Riemannian metric $(g_{ij})$ should be defined in a given tensor space $T^{0,2}$, and $H$ stands for the trace of $(h_{ij})$
\begin{equation}
H=g^{ij}h_{ij}.
\end{equation}
\ed

\br
 Functions $F$ that can be written as
\begin{equation}
F=GK^a,\q a>0,
\end{equation}
where $K$ is the Gaussian curvature and $G$ an arbitrary function of class $(K)$, including the case $G=1$, which doesn't belong to $(K)$, are of class $(K^*)$, \cf \cite[Proposition 1.9]{cg:indiana}.
\er

Thus, we shall solve the original problem \re{5.2} for curvature functions $F$ satisfying the requirement that their inverses $\tilde F\in (K^*)$.

Notice that in case $F=K$ there holds
\begin{equation}
F=\tilde F\in (K^*).
\end{equation}

We shall also assume without loss of generality that $F$ is homogeneous of degree $1$, and hence concave, \cf \cite[Lemma 1.2]{cg:minkowski}.

Now that we have formulated the condition for $\tilde F$, let us switch notations to enhance the readability of the text and to simplify the comparison with former results, and let us rewrite the equation \re{5.4} in the form
\begin{equation}\lae{5.9}
\fv FM=f(x)\qq\A\, x\in M,
\end{equation}
where $F\in (K^*)$, $0<f$ is defined in $N$  and $M\su N$ is a closed, spacelike, connected, strictly convex hypersurface, where its second fundamental form is defined with respect to the future directed normal, in contrast to our default convention to consider the past directed normal. 

In order to make the comparison with former results and techniques easier, we therefore switch the light cone, so that the future directed normal is now past directed, and replace the time function $\tau$ in \re{2.6} by $-\tau$ without changing the notation, i.e., $x^0$ is still the time function inherited from $\R[n+1,1]$, but 
\begin{equation}
\frac {d x^0}{d\tau}=-(1+\abs{x^0}^2),
\end{equation}
and the coordinate slices with positive curvature are now contained in $\{\tau<0\}$.

We want to solve equation \re{5.9}. For technical reasons, it is convenient to solve instead the
equivalent equation
\begin{equation}\lae{3.2}
\fv{\F(F)}M=\F(f),
\end{equation}
where $\F$ is a real function defined on $\R[]_+$ such that
\begin{equation}
\dot\F>0\q \tup{and}\q \ddot\F\le 0.
\end{equation}

For notational reasons, let us abbreviate
\begin{equation}
\tilde f=\F(f).
\end{equation}

We also point out that we may---and shall---assume without loss of generality
that $F$ is homogeneous of degree 1.

To solve \re{3.2} we look at the evolution problem
\begin{equation}\lae{3.5}
\begin{aligned}
\dot x&=(\F-\tilde f)\n,\\
x(0)&=x_0,
\end{aligned}
\end{equation}
where $x_0$ is an embedding of an initial strictly convex, compact, space-like
hypersurface $M_0$, $\F=\F(F)$, and $F$ is evaluated at the principal curvatures
of the flow hypersurfaces $M(t)$, or, equivalently, we may assume that $F$
depends on the second fundamental form $(h_{ij})$ and the metric $(g_{ij})$ of
$M(t)$; $x(t)$ is the embedding of $M(t)$, and $\nu$ is the past directed normal of the flow hypersurfaces $M(t)$.

This is a parabolic problem, so short-time existence is guaranteed---the proof
in the Lorentzian case is identical to that in the Riemannian case, cf. \ci[p.
622]{cg96}---, and under suitable assumptions, which we are going to formulate in a moment, we shall be able to prove that the
solution exists for all time and converges to a stationary solution if $t$ goes to
infinity.

In $N$ we consider an open, connected, precompact set $\Omega$ that is bounded by two
\tit{achronal}, connected, spacelike hypersurfaces $M_1\text{ and }M_2$, where
$M_1$ is supposed to lie in the \tit{past} of $M_2$.

We assume that $0<f\in C^{2,\al}(\bar \Om)$, $0<\al<1$, and that  the
boundary components $M_i$ act as barriers for $(F,f)$.

\bd\lad{0.1}
$M_2$ is an \tit{upper barrier} for $(F,f)$, if $M_2$ is strictly convex and
satisfies
\begin{equation}
\fm 2\ge f,
\end{equation}
and  $M_1$ is a \tit{lower barrier} for $(F,f)$, if at the points $\varSigma\su M_1$, where
$M_1$ is strictly convex, there holds
\begin{equation}
\fv F\varSigma\le f.
\end{equation}
$\varSigma$ may be empty.
\ed

To simplify some calculations that are to follow, we introduce an eigen time coordinate system in $N$, i.e., we write the metric in the form
\begin{equation}\lae{5.17c}
d\bar s^2=-d\tau^2+\cosh^2\tau \s_{ij}d\xi^id\xi^j,
\end{equation}
where $\s_{ij}$ is the standard metric of $S^n$. The time function $\tau$ is globally defined, and due to our convention the uniform convex slices are contained in $\{\tau<0\}$.

This preceding relation can be immediately deduced from \re{2.5} and \re{2.6}.  The special form of the metric with $\cosh^2\tau$ is of no importance. The crucial facts are that $N$ has constant curvature, $\pde{}\tau$ is a timelike unit vector field, and the coordinate slices $\{\tau=\const\}$ are totally umbilic.

Notice also that, if $M=\graph u$ is a spacelike hypersurface, the previously defined quantities $v$ and $\tilde v$ are identical to those defined in the new coordinate system
\begin{equation}
v^2=1-\bar g^{ij} u_iu_j,\q \bar g_{ij}= \cosh^2\tau\s_{ij}.
\end{equation}

However, when applying the formulas in \rs 1 one should observe that in the present coordinate system  the terms in equation \re{1.8b} should read 
\begin{equation}
\psi=0\q\wed\q \s_{ij}=\bar g_{ij}.
\end{equation}

We now consider the evolution problem \re{3.5} with $M_0=M_2$. Then the flow exists in a maximal time interval $I=[0,T^*)$, $0<T^*\le\un$, and, as we have proved in \cite[Section 4]{cg:indiana}, there holds
\bl
During the evolution the flow hypersurfaces stay inside $\bar \Om$. The hypersurfaces $M(t)$ can be written as graphs over $S^n$
\begin{equation}
M(t)=\graph u(t,\cdot),
\end{equation}
such that, if the barriers are expressed as $M_i=\graph u_i$, $i=1,2$, we have
\begin{equation}
u_1\le u\le u_2
\end{equation}
and the quantity
\begin{equation}
\tilde v=\frac1{\sqrt{1-\abs{Du}^2}}
\end{equation}
is uniformly bounded for all $t\in I$.

Moreover, the initial inequality $F\ge f$ is valid throughout the evolution, which can be equivalently formulated as
\begin{equation}\lae{5.23c}
\F\ge \tilde f.
\end{equation}
\el

Let us now look at the evolution equations satisfied by $u$, $\tilde v$, and $h^j_i$.

\bl
Let $M(t)=\fv{\graph u(t,\cdot)}{S^n}$ be the flow hypersurfaces, then $u$ satisfies the parabolic equation
\begin{equation}\lae{5.24c}
\dot u-\dot\F F^{ij}u_{ij}=-\tilde v (\F-\tilde f)+\tilde v\dot\F F -\dot\F F^{ij}\bar h_{ij},
\end{equation}
where
\begin{equation}
\dot\F=\frac{d\F(r)}{dr}\q\wed\q \dot u=\tfrac{du}{dt}.
\end{equation}
\el

\bp
These equations immediately follow from the relations \re{1.16}, \re{1.18} and \re{3.5} by observing that
\begin{equation}
\dot u=\dot x^0.
\end{equation}
Notice that $\dot u$ is the total time derivative, where \cq{time} is just the usual name for the flow parameter.
\ep

\bl
The quantity $\tilde v$ satisfies the evolution equation
\begin{equation}\lae{5.24}
\begin{aligned}
\dot{\tilde v}-\dot\F F^{ij}\tilde v_{ij}&=  -\dot\F F^{ij}h_{ik}h^k_j\tilde v -[(\F-\tilde f)-\dot\F F] \bar h_{ij}\check u^i\check u^j\tilde v^2\\
&\hp{=}\;  +2 \dot \F F^{ij}h^k_i\bar h_{ki}
-\bar\ka^2\dot\F F^{ij}\bar g_{ij}\tilde v  -2 \bar\ka^2\dot\F F^{ij}u_iu_j\tilde v \\
&\hp{=}\;+K_N\dot\F F^{ij}u_iu_j\tilde v+f_\bet x^\bet_iu^i,
\end{aligned}
\end{equation}
where $\bar \ka$ is the principal curvature of the slices $\{\tau=\const\}$, 
\begin{equation}
\check u^i=\bar g^{ij}u_j 
\end{equation}
and $K_N=1$ for the de Sitter space-time.
\el

\bp
Let $(\h_\al)=(-1,0,\ldots,0)$ be the covariant vector field representing $-\pde{}\tau$. Differentiating $\tilde v=\h_\al\nu^\al$ covariantly with respect to the induced metric of $M$, where $(\nu^\al)$ is the past directed normal, we obtain
\begin{equation}\lae{5.29c}
\begin{aligned}
\dot{\tilde v}-\dot\F F^{ij}\tilde v_{ij}=&-\dot\F F^{ij}h_{ik}h_j^k\tilde v
+[(\F-\tilde f)-\dot\F F]\h_{\alpha\beta}\n^\alpha\n^\beta\\
&-2\dot\F F^{ij}h_j^k x_i^\alpha x_k^\beta \h_{\alpha\beta}-\dot\F F^{ij}\h_{\alpha\beta\gamma}x_i^\beta
x_j^\gamma\n^\alpha\\
&-\dot\F F^{ij}\riema \alpha\beta\gamma\delta\n^\alpha x_i^\beta x_k^\gamma x_j^\delta\h_\e x_l^\e g^{kl}\\
&-\tilde f_\beta x_i^\beta x_k^\alpha \h_\alpha g^{ik},
\end{aligned}
\end{equation}
where the ambient space can be a general Lorentzian manifold, \cf \cite[Lemma 4.4]{cg:indiana}; however, in the general case the definition of $\h$ has to be adjusted, since it has to be a unit vector field.

Now, if the ambient space is a space of constant curvature $K_N$, the term containing the Riemannian curvature tensor vanishes, and we shall show that the crucial term
\begin{equation}
-F^{ij}\h_{\al\bet\ga}x^\bet_ix^\ga_j\nu^\al
\end{equation}
can be expressed as claimed.

First we observe that the second fundamental form $\bar h_{ij}$ of the coordinate slices $\{\tau=\const\}$ is given by
\begin{equation}\lae{5.28}
\begin{aligned}
\bar h_{ij}&=-\tfrac12\dot {\bar g}_{ij}.
\end{aligned}
\end{equation}

Secondly, the vector field $(\h_\al)$ is a gradient field, namely,
\begin{equation}
(\h_\al)=\grad\f
\end{equation}
with
\begin{equation}
\f(\tau,x^i)=-\tau.
\end{equation}

Since $D\f$ is a unit vector field, we have
\begin{equation}\lae{5.32}
\f_{\al\bet}\f^\al=\f_{\bet\al}\f^\al=0.
\end{equation}

The restriction of $\f$ to a coordinate slice is constant, hence, differentiating $\f$ covariantly with respect to the induced metric $\bar g_{ij}$, we deduce

\begin{equation}
\begin{aligned}
0=\f_{ij}&=\f_\al\bar x^\al_{ij}+\f_{\al\bet}\bar x^\al_i\bar x^\bet_j\\
&=\f_\al\bar\nu^\al\bar h_{ij}+\f_{\al\bet}\bar x^\al_i\bar x^\bet_j,
\end{aligned}
\end{equation}
where
\begin{equation}
\bar x=(\tau,x^i,\ldots,x^n),\q\tau=\const,
\end{equation}
is the embedding of the coordinate slice, and we conclude
\begin{equation}\lae{5.35}
\f_{\al\bet}\bar x^\al_i\bar x^\bet_j=-\f_\al\bar\nu^\al\bar h_{ij}=-\bar h_{ij}
\end{equation}
as well as
\begin{equation}\lae{5.36}
\f_{\al\bet}\bar\nu^\bet=0.
\end{equation}

Differentiating \re{5.35} covariantly with respect to the induced metric of the coordinate slices, we infer
\begin{equation}\lae{5.37}
\f_{\al\bet\ga}\bar x^\al_i\bar x^\bet_j\bar x^\ga_k=-\bar h_{ij;k}=0,
\end{equation}
in view of \re{5.28}, where we also used \re{5.36}.

Finally, differentiating \re{5.32} covariantly, we conclude
\begin{equation}\lae{5.38}
0=\f_{\al\bet\ga}\f^\al+\f_{\al\bet}\f^\al_\ga.
\end{equation}

We can now evaluate the term
\begin{equation}
-\h_{\al\bet\ga}\nu^\al x^\bet_ix^\ga_j=-\f_{\al\bet\ga}\nu^\al x^\bet_ix^\ga_j
\end{equation}
with the help of the relations \re{5.37} and \re{5.38} yielding
\begin{equation}
\begin{aligned}
-\f_{\al\bet\ga}\nu^\al x^\bet_ix^\ga_j&=-\f_{0ij}\nu^0+\tilde v\f_{k0j}\check u^ku_i+\tilde v\f_{ki0}\check u^ku_j\\
&=-\tilde v \bar h^k_i\bar h_{kj}-\tilde v  \bar h^m_k\bar h_{mi}\check u^ku_j-\tilde v \bar h^m_k\bar h_{mj}\check u^ku_i\\
&\hp{=}\;+K_N\tilde v u_iu_j,
\end{aligned}
\end{equation}
where we applied the Ricci identities.

The indices of $\bar h_{ij}$ are raised with the help of the metric $\bar g_{ij}$.

Taking then \re{5.28} into account completes the proof of the lemma. 
\ep

The equation \re{5.24} can be even simplified further by using the same argument as in the case of its Riemannian analogue, \cf \cite[Lemma 5.8]{cg:spaceform}.

Let $\h=\h(\tau)$ be a positive solution of the ordinary differential equation
\begin{equation}
\dot\h=-\tfrac{\bar H}n\h=-\bar\ka \h,
\end{equation}
notice that $\h$ is defined for all $\tau\in\R[]$, and set 
\begin{equation}\lae{5.44c}
\chi=\tilde v\h.
\end{equation}

Then we can prove

\bl
The function $\chi$ satisfies the evolution equation
\begin{equation}\lae{5.45c}
\begin{aligned}
\dot\chi -\dot\F F^{ij}\chi_{ij}=-\dot\F F^{ij}h_{ki}h^k_j\chi +[(\F-\tilde f)+\dot\F F]v \bar \ka \chi+f_\al x^\al_iu^iv\chi,
\end{aligned}
\end{equation}
for any value of $K_N$.
\el
\bp
Differentiating \re{5.44c} we deduce
\begin{equation}\lae{5.46c}
\begin{aligned}
\dot\chi-\dot\F F^{ij}\chi_{ij}&=\{\dot {\tilde v} -\dot\F F^{ij}\tilde v_{ij}\}\h+\{\dot u -\dot\F F^{ij}u_{ij}\}\tilde v \dot \h\\
&\hp{=}\;-2\dot\h \dot \F F^{ij}\tilde v_i u_j -\tilde v\Ddot\h\dot\F F^{ij} u_iu_j.
\end{aligned}
\end{equation}

Now we first observe
\begin{equation}\lae{5.47c}
\begin{aligned}
\Ddot\h&=-\frac{\dot{\bar H}}n\h-\frac{\bar H}n\dot\h=-\tfrac1n\{\abs{\bar A}^2+\bar R_{\al\bet}\bar\nu^\al\bar\nu^\bet\}\h+\bar\ka^2\h\\
&=K_N\h.
\end{aligned}
\end{equation}

Secondly, from 
\begin{equation}
\tilde v^2 =1+\norm{Du}^2
\end{equation}
we derive
\begin{equation}
\begin{aligned}
\tilde v_i= v u_{ij}u^j&= -h_{ij}u^j+\bar h_{ij}u^j v\\
&=-h_{ij}u^j+\bar\ka u_i\tilde v,
\end{aligned}
\end{equation}
hence we obtain
\begin{equation}\lae{5.50c}
\begin{aligned}
-2\dot\h\dot\F F^{ij}\tilde v_iu_j&=-2\bar\ka \h \dot\F F^{ij}h_{ik}u^ku_j+2\bar\ka^2\dot F^{ij}u_iu_j\tilde v\h.
\end{aligned}
\end{equation}

Inserting \re{5.47c} and \re{5.50c} in \re{5.46c} we conclude
\begin{equation}
\begin{aligned}
\dot\chi-\dot\F F^{ij}\chi_{ij}&=\{\dot {\tilde v} -\dot\F F^{ij}\tilde v_{ij}\}\h-\bar \ka\{\dot u -\dot\F F^{ij}u_{ij}\}\tilde v \h\\
&\hp{=}\;-2\bar\ka\dot  \F F^{ij}h_{ik}u^ku_j\h +2\bar\ka^2\dot\F F^{ij}u_iu_j\tilde v\h\\
&\hp{=}\;-K_N\dot\F F^{ij}u_iu_j\tilde v\h,
\end{aligned}
\end{equation}
from which the result immediately follows, in view of \re{5.24c} and \re{5.24}.
\ep

\br\lar{5.8}
Since the flow stays in a compact subset and the hypersurfaces $M(t)$ are uniformly convex, there exist positive constants $c_1$, $c_2$ such that
\begin{equation}
0<c_1\le \chi\le c_2.
\end{equation}
 This follows immediately from the observation that in a point, where $D\chi=0$, there holds
 \begin{equation}
h_{ij}u^j=0,
\end{equation}
 hence $Du=0$, and thus $\chi=\h$.
\er
In case of a Lorentzian space form the evolution equation for the second fundamental form  is a rather simple expression.
\bl\lal{5.9}
The second fundamental form $(h^j_i)$ satisfies the differential equation
\begin{equation}\lae{5.42}
\begin{aligned}
\dot h_i^j-\dot\F F^{kl}h_{i;kl}^j&=-\dot\F F^{kl}h_{rk}h_l^rh_i^j+\dot\F F
h_{ri}h^{rj}- (\F-\tilde f) h_i^kh_k^j\\
&\hp{=}\;-\tilde f_{\alpha\beta} x_i^\alpha x_k^\beta g^{kj}- \tilde f_\alpha\n^\alpha h_i^j+\dot\F
F^{kl,rs}h_{kl;i}h_{rs;}^{\hphantom{rs;}j}\\
&\hp{=}\;+\Ddot\F F_iF^j\\
&\hp{=}\;+K_N\{(\F-\tilde f)+\dot\F F\de^j_i-\dot\F F^{kl}g_{kl}h^j_i\}.
\end{aligned}
\end{equation}
\el

\bp
The evolution equation for $h^j_i$ in a semi-Riemannian manifold has been derived in \cite[Lemma 3.5]{cg:indiana}.

If the ambient space is a Lorentzian space form  $N$ with curvature $K_N$, and if $M$ is a spacelike hypersurface with normal $\nu$,  then the former  equation reduces to \re{5.42}.

For Riemannian space forms,  this equation has already been established in \cite[Corollary 5.4]{cg:spaceform}. Of course, in the stationary case, this identity was first proved by J. Simons in \cite{simons}.
\ep

\section{Curvature estimates}\las 6
We are now able to prove the a priori estimate for the principal curvatures of the $M(t)$.

\bl\lal{6.1c}
Consider the flow in a maximal interval $I=[0,T^*)$, choose $\F=\log$,  and assume that the initial hypersurface  $M_2$ is of of class $C^{4,\al}$, where $F\in (K^*)$ and $0<f\in C^{2,\al}(\bar \Om)$. Then there are positive constants $k_1,k_2$, depending only on $F$, $f$ and $\Om$, such that the principal curvatures $\ka_i$ are estimated by
\begin{equation}
0<k_1\le\ka_i\le k_2.
\end{equation}
\el

\bp
It suffices to prove an upper estimate for $\ka_i$, since $\fv F{\pa\C_+}=0$.

We observe that $u$, $\tilde v$ and $\chi$ are already uniformly bounded, and that $\chi$ is also uniformly positive, \cf \rr{5.8}.

Let $\f$ and $w$ be defined respectively by
\begin{align}
\f&=\sup\set{{h_{ij}\h^i\h^j}}{{\norm\h=1}},\\
w&=\log\f+\lambda\chi,\lae{6.3c}
\end{align}
where $\lambda$ is a large positive parameter to be specified later. We claim that
$w$ is bounded for a suitable choice of $\lambda$.

Let $0<T<T^*$, and $x_0=x_0(t_0)$, with $ 0<t_0\le T$, be a point in $M(t_0)$ such
that

\begin{equation}
\sup_{M_0}w<\sup\set {\sup_{M(t)} w}{0<t\le T}=w(x_0).
\end{equation}

We then introduce a Riemannian normal coordinate system $(\x^i)$ at $x_0\in
M(t_0)$ such that at $x_0=x(t_0,\x_0)$ we have
\begin{equation}
g_{ij}=\delta_{ij}\q \tup{and}\q \f=h_n^n.
\end{equation}

Let $\tilde \h=(\tilde \h^i)$ be the contravariant vector field defined by
\begin{equation}
\tilde \h=(0,\dotsc,0,1),
\end{equation}
and set
\begin{equation}
\tilde \f=\frac{h_{ij}\tilde \h^i\tilde \h^j}{g_{ij}\tilde \h^i\tilde \h^j}\raise 2pt
\hbox{.}
\end{equation}

$\tilde \f$ is well defined in neighbourhood of $(t_0,\x_0)$.

Now, define $\tilde w$ by replacing $\f$ by $\tilde \f$ in \re{6.3c}; then, $\tilde w$
assumes its maximum at $(t_0,\x_0)$. Moreover, at $(t_0,\x_0)$ we have
\begin{equation}
\dot{\tilde \f}=\dot h_n^n,
\end{equation}
and the spatial derivatives do also coincide; in short, at $(t_0,\x_0)$ $\tilde \f$
satisfies the same differential equation \re{5.42} as $h_n^n$. For the sake of
greater clarity, let us therefore treat $h_n^n$ like a scalar and pretend that $w$
is defined by
\begin{equation}
w=\log h_n^n+\lambda\chi.
\end{equation}

At $(t_0,\x_0)$ we have $\dot w\ge 0$, and, in view of the maximum principle,  we deduce from \re{5.5}, \re{5.42}, and \re{5.45c}
\begin{equation}\lae{6.10c}
\begin{aligned}
0\le&\msp[3] \dot\F F h_n^n-(\F-\tilde f) h_n^n+\lambda c_1 -\lambda\e_0\dot\F F H \chi\\ & +\lambda c_1[(\F-\tilde f)+\dot\F F]\\
&+\dot\F F^{ij}(\log h_n^n)_i(\log h_n^n)_j\\
&+\{\ddot\F F_n F^n +\dot\F
F^{kl,rs}h_{kl;n}h_{rs;}^{\hphantom{rs;}n}\}(h_n^n)^{-1},
\end{aligned}
\end{equation}
where we have estimated bounded terms by a constant $c_1$, assumed that
$h_n^n, \lambda$ are larger than $1$, and used \re{5.5} as well as the simple
observation
\begin{equation}\lae{5.14}
\abs{F^{ij}h_j^k\h_{ik}}\le \norm\h F
\end{equation}
valid for any tensor field $(\h_{ik})$.

Now, the last term in \re{6.10c} is estimated from above by
\begin{equation}\lae{6.12c}
\{\ddot\F F_n F^n+\dot\F F^{-1} F_n F^n\}(h_n^n)^{-1}-\dot \F F^{ij}
h_{in;n}h_{jn;}^{\hphantom{jn;}n}(h_n^n)^{-2},
\end{equation}
\cf \cite[Lemma 1.5]{cg:indiana}, where  the sum in the braces vanishes,
due to the choice of
$\F$. Moreover, because of the Codazzi equation, we have
\begin{equation}
h_{in;n}=h_{nn;i},
\end{equation}
and hence,  we conclude
that \re{6.12c} is bounded from above by
\begin{equation}
-(h_n^n)^{-2}\dot\F F^{ij}h_{n;i}^nh_{n;j}^n.
\end{equation}

Thus, the terms in \re{6.10c} containing the derivatives of $h_n^n$ sum up to something non-positive.

Choosing then in \re{6.10c} $\lambda$ such that
\begin{equation}
2\le \lambda\e_0\chi
\end{equation}
we derive
\begin{equation}
\begin{aligned}
0\le&\msp[2] -\dot\F F H -(\F-\tilde f) h_n^n\\
&+\lambda c_1 [(\F-\tilde f)+\dot\F F]+\lambda c_1.
\end{aligned}
\end{equation}

We now observe that $\dot \F F=1$, and deduce in view of \re{5.23c} that $h_n^n$ is
a priori bounded at $(t_0,\x_0)$.
\ep

The result of the preceding lemma can be restated as a uniform estimate for the functions
$u(t)\in C^2(S^n)$. Since, moreover, the principal curvatures of the flow
hypersurfaces are not only bounded, but also uniformly bounded away from zero,
in view of \re{5.23c} and the assumption that $F$ vanishes on $\pa \C_+$, we
conclude that $F$ is uniformly elliptic on $M(t)$.

\section{Convergence to a stationary solution}\las 7

We are now ready to prove \rt{0.4}. Let $M(t)$ be the flow with initial
hypersurface $M_0=M_2$. Let us look at the scalar version of the flow 
\begin{equation}\lae{7.1}
\pde ut=-v(\F-\tilde f),
\end{equation}
\cf \cite[equ. (3.5)]{cg:indiana}. 

This is  a scalar parabolic differential equation defined on the cylinder
\begin{equation}
Q_{T^*}=[0,T^*)\times S^n
\end{equation}
with initial value $u(0)=u_2\in C^{4,\alpha}(S^n)$. In view of the a priori estimates,
which we have established in the preceding sections, we know that
\begin{equation}
{\abs u}_\low{2,0,S^n}\le c
\end{equation}
and
\begin{equation}
\F(F)\,\tup{is uniformly elliptic in}\,u
\end{equation}
independent of $t$. Moreover, $\F(F)$ is concave, and thus, we can apply
the regularity results of \ci[Chapter 5.5]{nk} to conclude that uniform
$C^{2,\alpha}$-estimates are valid, leading further to uniform $C^{4,\alpha}$-estimates
due to the regularity results for linear operators.

Therefore, the maximal time interval is unbounded, i.e. $T^*=\un$.

Now, integrating \re{7.1} with respect to $t$, and observing that the right-hand
side is non-positive, yields
\begin{equation}
u(0,x)-u(t,x)=\int_0^tv(\F-\tilde f)\ge c\int_0^t(\F-\tilde f),
\end{equation}
i.e.,
\begin{equation}
\int_0^\un \abs{\F-\tilde f}<\un\qq\A\msp x\in S^n
\end{equation}
Hence, for any $x\in S^n$ there is a sequence $t_k\rightarrow \un$ such that
$(\F-\tilde f)\rightarrow 0$.

On the other hand, $u(\cdot,x)$ is monotone decreasing and therefore
\begin{equation}
\lim_{t\rightarrow \un}u(t,x)=\tilde u(x)
\end{equation}
exists and is of class $C^{4,\alpha}(S^n)$ in view of the a priori estimates. We, finally,
conclude that $\tilde u$ is a stationary solution of our problem, and that
\begin{equation}
\lim_{t\rightarrow \un}(\F-\tilde f)=0.
\end{equation}

%\backmatter

\bibliographystyle{hamsplain}
%\bibliography{mrabbrev,publications}
 \providecommand{\bysame}{\leavevmode\hbox to3em{\hrulefill}\thinspace}
\providecommand{\href}[2]{#2}

%\listoffigures

%\cleardoublepage

%\thispagestyle{empty}
%\closegraphsfile
\end{document}